\documentclass[preprint,12pt]{elsarticle}
\usepackage{amsmath,latexsym,amssymb, enumerate, amsthm, epsfig }
\usepackage{amsmath,amsfonts,amsthm,amscd,amsxtra,enumerate}
\usepackage[paperheight=297mm, paperwidth=210mm, left=30mm, right=20mm]{geometry}
\usepackage{hyperref}
\usepackage[hyperref]{xcolor}
\hypersetup{colorlinks=true,}

\begin{document}

\begin{frontmatter}

\title{Analytical and numerical study of weakly nonlinear hyperbolic waves in a van der Waals gas}


\author{Harsh V. Mahara and V. D. Sharma}

\address{Department of Mathematics, Indian Institute of Technology Bombay,
	Powai, Mumbai-400076}

\begin{abstract}
In this paper, we characterized resonant interaction of weakly nonlinear
hyperbolic waves in gas dynamics with a real gas background. An asymptotic approach
is used to study the interaction between waves, governed by the Euler equations of gas
dynamics, supplemented by a van der Waal equation of state; the evolution equation is an
integro-differential equation composed of a Burgers type nonlinear term and a convolution
term with a known kernel. A one parameter family of traveling wave solutions for
various values of van der Waals parameter is studied analytically and numerically. Effect
of the influence of van der Waals parameter on the properties of traveling waves solution
is investigated. Numerical experiments on the evolution of arbitrary initial data are
performed using fractional step algorithm, describing the behavior of convolution term
in the integrodiffrential equation. The existence of non breaking for all time solutions is
substantiated numerically.

\end{abstract}


\end{frontmatter}


\section{Introduction}
In one-dimensional case, the study of the initial stages of the hyperbolic wave phenomena with properly confined initial data is well understood ( Courant and Friedrichs \cite{MR0029615},
DiPerna \cite{MR338576}, Glimm \cite{MR0194770}, Lax \cite{Lax} , Nishida and Smoller \cite{MR0330789}, and Whitham \cite{GB74}). In an open domain, after some time the component waves of the solution (moving with different speed) get separated in the space and hence the nonlinear interaction between them comes to an end ( Liu \cite{MR0450802}). On the contrary, the understanding of the long-time nature of the solution to the compressible fluid flow in bounded regions (periodic case) is still incomplete.

Majda and Rosales \cite{MR760229} first developed the perturbation theory for resonantly interacting weakly nonlinear hyperbolic waves and reduced the system of one-dimensional system of Euler equations of gasdynamics into a pair of integro-differential equation having inviscid Burger's equation as nonlinear component coupled through a linear integral operator with the known kernel. In the numerical experiment, Majda, Rosales and Schnobek \cite{MR975485} obtained a traveling wave solution and a wave train eliminating completely the shock front. Pego \cite{MR975486} and Gabov \cite{MR510246} found an exact expression for standing traveling wave solution. The surprising existence of the other wave is explained as the interaction of the acoustic components through the entropy variation might balance the nonlinear effects and gives rise to the continuous nontrivial waveform. 

For Euler equations, the existence of non breaking (no shock) for all time solution has long been an open problem. Motivated by these results Celentano \cite{ccm}, found a large class of periodic standing waves using bifurcation theory and Vaynblat \cite{VD96} showed that there are other complex structured new never breaking solutions exists for complete gas dynamics equations. Shefter and Rosales \cite{MR1719749} considered in the ideal gas background the one-dimensional Euler equations of gas inside a closed duct in the weakly nonlinear regime. They showed analytical and numerical evidence for the existence of new quasiperiodic and globally attracting nonbreaking for all time (NBAT) solution for the one-dimensional inviscid Euler equations of gasdynamics.

In this chapter, we extend their work to the real gas regime and using asymptotic analysis reduce the one-dimensional  Euler's equation of gas dynamics with periodic boundary condition to a single leading order asymptotic integro-differential  evolution equation 

\begin{equation}
\frac{\partial}{\partial t} \sigma(x,t)+
\Lambda \frac{\partial}{\partial x}\left[\frac{1}{2}\sigma^2(x,t)\right] + \Gamma
\int_{0}^{2\pi}K(x-y)
\sigma(y,t)~dy=0,
\end{equation}

where $ \sigma (x,t) $ is the acoustic amplitude, and $ K(x) $ is the kernel of the convolution term and is equal to the derivative of the entropy and the coefficients $ \Lambda $ and $ \Gamma $ correspond to the real gas effects and are defined later in \eqref{a.2.20}. In all the numerical experiments we consider a single-mode kernel $ K(x)=\sin x.$   

We have then analyzed  analytically and numerically the traveling wave solutions. The influence of van der Waal parameter $ b $ on these solutions has been studied. We have obtained  non breaking for all time solutions for real gases and have shown their globally attracting behavior by numerical experiments. For our numerical experiments
we have used second-order accurate fractional step approach \cite{MR0235754}.

The work is organized as follows: Basic equations and formulation of the problem are given in Section \ref{basic equation}. In Section \ref{evolution equation}, a detailed derivation of a pair of evolution equations is given and a precise expression for a family of traveling wave solutions for this system  is obtained. Then this pair of equations is reduced to a single evolution equation which is the main object of study all over the chapter. We study the influence of real gas effects on the traveling wave solution in Section \ref{trav}. In Section \ref{num}, we give a brief review of the numerical technique used in the experiment. Evolution of arbitrary initial data for different values of  van der Waal parameter is investigated the presence of the NBAT solution is noticed. Finally, we conclude this chapter with a discussion of our results in Section \ref{con}.

\section{Basic equations}\label{basic equation}
In this section, we study the motion of a quiescent real gas inside a rigid pipe with close ends and having, a circular cross-section of constant radius. The basic equations governing the unsteady, inviscid, one-dimensional flow of a real fluid obeying van der Waal equation of state are the Euler equations of gasdynamics which are given by 
\begin{equation}\label{a.2.1}
\begin{aligned}
&\bar{\rho}_{\bar{t}}+(\bar{\rho}\bar{v})_{\bar{x}}=0,\\
&(\bar{\rho}\bar{v})_{\bar{t}}+(\bar{\rho}\bar{v}^2+\bar{P})_{\bar{x}}=0,\\
&(\bar{\rho}[\bar{e}+\frac{\bar{v}^{2}}{2}])_{\bar{t}}+
(\bar{\rho}\bar{v}[\bar{e}+\frac{\bar{P}}{\bar{\rho}}+\frac{\bar{v}^{2}}{2}])_{\bar{x}}=0,
\end{aligned}	
\end{equation}	
where  $ \bar{\rho} $ is the mass density, $ \bar{v} $ is the flow velocity, $ \bar{P} $ is the gas pressure, and $ \bar{e} $ is internal energy per unit mass. The variables $ \bar{x} $ and $ \bar{t} $ are space coordinate and time, respectively; the subscripts denote the partial derivatives.	

Equations \eqref{a.2.1} are supplemented by an equation of state which in our case is van der Waals equation of state of the form \cite{MR1419777}

\begin{equation}\label{a.2.2}
\bar{P} = K_{0}\delta\frac{\bar{\rho}^{1+\delta}\exp(\delta\bar{s}/R)}{(1-\bar{b}\bar{\rho})^{1+\delta}},
\qquad
\bar{T} = K_{0}\delta\frac{\bar{\rho}^{\delta}\exp(\delta\bar{s}/R)}{R(1-\bar{b}\bar{\rho})^{\delta}},
\end{equation}
where, $ K_{0} $ is a constant, $ \bar{s} $ is the entropy, $ \bar{T} $ is the temperature, $ \delta $ is a dimensionless material-dependent quantity defined as $ \delta = R/c_{v} $ with $ c_{v} $ the specific heat at constant volume and $ R $ the specific gas constant;
$ \delta $ lies in the interval $ 0 < \delta \leq 2/3  $ with $  \delta= 2/3  $ for a monoatomic fluid and the parameter $ \bar{b} $ represents the van der Waal excluded volume.

In view of the thermodynamic identity
\begin{equation}\label{a.2.3}
\displaystyle{\bar{T} d\bar{s} = d \bar{e} + \bar{P} d \left(\frac{1}{\bar{\rho}}\right)},
\end{equation} 
satisfied by the gas variables, we obtain the constitutive relations 
$ \bar{P}=\bar{\rho}^{2}\bar{e}_{\bar{\rho}},\; \bar{P}=\bar{e}_{\bar{s}} $ and
$ \bar{P}_{\bar{s}}=\bar{\rho}^2 \bar{T}_{\bar{\rho}} $. Further, the sound speed is defined by $ \bar{c}=(\bar{P}_{\bar{\rho}})^{1/2}>0.$
In addition, we assume that the gas is at rest at the closed ends, i.e., the boundary conditions are given by
$$ \bar{v}(0,\bar{t})=0  ~\text{ and }~  \bar{v}(l,\bar{t})=0  .$$

We introduce non-dimensional variables, defined as
\begin{equation*}
x=\frac{\pi}{l}\bar{x},\;\; t=\frac{\pi\sqrt{\bar{P}_{0}/ \bar{\rho}_{0}}}{l}\bar{t},
\;\; P=\frac{\bar{P}}{\bar{P}_{0}}, \;\; \rho=\frac{\bar{\rho}}{\bar{\rho}_{0}},
\;\; T=\frac{\bar{T}}{\bar{T}_{0}}, \;\; s=\frac{\bar{s}}{\bar{s}_{0}},
\;\; v=\frac{\bar{v}}{\sqrt{\bar{P}_{0}/ \bar{\rho}_{0}}},
\;\; b = \bar{b}\bar{\rho}_{0},
\end{equation*}
where $  l $ is the chosen length of the tube. Using  dimensionless variables  in \eqref{a.2.1} and in other previous thermodynamic relations, we see that they
remain unchanged in the dimensionless form. We choose $ \rho_{0} $ so that the mean value of $ \rho  $ is normalized to 1 and $ s=0 $;  finally, we choose $ P_{0}~\text{and}~ T_{0} $ so that $ c= \sqrt{(1 + \delta )/ (1-b) } $ and replace
the tube length $ l $ by $ \pi $ in the boundary conditions.

The problem of a gas inside a closed tube is equivalent to a spatially periodic problem with a period twice the length of the tube. For every solution of the original problem, we can extend it to a solution of the equation defined for all $ x $ periodic with period $ 2\pi $; with $ P $ and $ \rho $ ( thus $ s $, $ T $, and $ e $) even and with $ v $ odd. Conversely, if one has a solution with these properties then it reduces to the solution of the tube problem above when restricted to the tube. Thus we substitute the boundary condition by analogous conditions:

\begin{equation}\label{a.2.4}
\begin{aligned}
&P(x+2\pi,t) = P(x,t), \quad\qquad \!\!\!\!P(x,t) = P(-x,t),\\
&\rho(x+2\pi,t) = \rho (x,t),\quad \qquad \rho(x,t) = \rho(-x,t),\\
&e(x+2\pi,t) = e (x,t), \quad\qquad e(x,t) = e (-x,t), \\
&s(x+2\pi,t) = s (x,t), \quad\qquad s(x,t) = s (-x,t), \\
&v(x+2\pi,t) = v (x,t), \quad\qquad \!v(x,t) = -v (-x,t), \\
\end{aligned}
\end{equation}
for every $ x $ and $ t, $ here any two rows of the first four rows imply the other two.

We introduce the state vector $ \mathbf{U} = (\rho, v, s)^{tr} $ so that \eqref{a.2.1} may be written in the vector matrix notation
\begin{align}\label{a.2.5}
\mathbf{U}_{t}+\mathbf{A}(\mathbf{U})\mathbf{U}_{x}=\mathbf{0},
\end{align}
where $ \mathbf{A} $ is a square matrix with components denoted by $  A_{ij}  $, having nonzero components
$$ A_{11}=A_{22}=A_{33}=v, \quad  A_{12}=\rho,\quad A_{21}=\frac{P_{\rho}}{\rho},\quad A_{21}=\frac{P_{s}}{\rho}. $$
We analyse the interaction of waves, which propagate through a constant background state
$ \mathbf{U}_{0} = (1, 0,0)^{tr} .$
The system \eqref{a.2.5} is strictly hyperbolic with characteristic velocities 
$\; \lambda_{1}=-c_{0},\;\lambda_{2}=0,\;\lambda_{3}=c_{0} ,$ and corresponding right eigenvectors are given by
\begin{equation}\label{l1.3}			
\mathbf{R}_{1}=\left(\begin{array}{ccccc}
1\\
-c_{0}\\
0 \\
\end{array}\right),\qquad
\mathbf{R}_{2}=\left(\begin{array}{ccccc}
P_{s0}\\
0\\
-c_{0}^2\\
\end{array}\right),
\qquad
\mathbf{R}_{3}=\left(\begin{array}{ccccc}
1\\
c_{0}\\
0\\
\end{array}\right).
\end{equation}
The associated left eigenvectors $ \mathbf{L}_{i},\; i = 1, 2, 3, $ can be obtained using the normalization condition
$  \mathbf{L}_{i}. \mathbf{R}_{i} = \delta_{ij}$ where $ \delta_{ij}  $ is kronecker delta.
Thus, the $ 1 $-wave and $ 3 $-wave are the left and right moving sound waves and the $ 2 $-wave corresponds to the convection of the entropy with the particle velocity. The subscript $ 0 $ refers to the evaluation at $ \mathbf{U} = \mathbf{U}_{0} $ and is synonymous with the state of equilibrium.
\section{Derivation of evolution equations}\label{evolution equation}
In this section, we give a derivation of a pair of evolution equations exhibiting the interaction between the acoustic components. We look for a uniformly valid asymptotic solution to \eqref{a.2.5} as $ \epsilon \to 0 $ of the form 
\begin{equation}\label{a.2.7}
\mathbf{U}= \mathbf{U}_{0}+\epsilon \mathbf{U}_{1}(x,t,\xi,\tau)+
\epsilon^2 \mathbf{U}_{2}(x,t,\xi,\tau)+O(\epsilon^3),
\end{equation}
where $ \displaystyle {\xi=\frac{x}{\epsilon},~ \tau=\frac{t}{\epsilon}} $ are fast variables.
The cofficient matrix $ \mathbf{A} $ can be extended about the constant state 
$ \mathbf{U}_{0} $
\begin{equation}\label{a.2.8}
\mathbf{A}(\mathbf{U})= \mathbf{A}(\mathbf{U}_{0})+
\nabla\mathbf{A}(\mathbf{U}_{0})(\mathbf{U}-\mathbf{U}_{0})+
\frac{1}{2}\nabla^2\mathbf{A}(\mathbf{U}_{0}).(\mathbf{U}-\mathbf{U}_{0})(\mathbf{U}-\mathbf{U}_{0})^{tr}+\ldots
\end{equation}
Using  \eqref{a.2.7} in \eqref{a.2.8}, and equating the coefficients of $ \epsilon $ we find that $ \mathbf{U}_{1}, \mathbf{U}_{2} $ satisfy
\begin{equation}\label{a.2.9}
\mathbf{U}_{1\tau}+\mathbf{A}_{0}\mathbf{U}_{1\xi}=0,
\end{equation}
\begin{equation}\label{a.2.10}
\mathbf{U}_{2\tau}+\mathbf{A}_{0}\mathbf{U}_{2\xi} +\mathbf{B}(x,t,\xi,\tau) = 0,
\end{equation}
where, $ \mathbf{A}_{0}=\mathbf{A}(\mathbf{U}_{0}),
\;\nabla\mathbf{A}_{0}=\nabla\mathbf{A}(\mathbf{U}_{0}),
~\text{and}~ \;\mathbf{B}= \mathbf{U}_{1t}+\mathbf{A}_{0}\mathbf{U}_{1x}+
\nabla\mathbf{A}_{0}\mathbf{U}_{1}\mathbf{U}_{1\xi}.$

A solution of \eqref{a.2.9} is given by 
\begin{equation}\label{a.2.11}
\mathbf{U}_{1}=\sum_{j=1}^{3}a_{j}(x,t,\theta_{j})\,\mathbf{R}_{j},
\end{equation}
where $ a_{j} = (\mathbf{L}_{j}.\mathbf{U}_{1}) $ is a scalar function called the wave amplitude that depends on the $ j $-th phase variable $ \theta_{j} $ defined as
$ \theta_{j} = k_{j}\xi -\omega_{j}\tau.$
The dependence of $ a_{j} $ on $ \theta $ describes the waveform; the wave number $ k_{j} $ and frequency $ \omega_{j} $ satisfy $\; \omega_{j}=\lambda_{j}k_{j},\; j=1,2,3. $
We suppose that the wave amplitudes and their derivatives with respect to $ x,\, t, \,\text{and}~ \theta $ are periodic or almost periodic functions of $ \theta .$ 
For simplicity, we also suppose that $ a_{j}(x,t,\theta) $ has zero mean with respect to the phase variable $ \theta_{j} ,$ i.e., 
\begin{equation*}
\lim_{T\to\infty}\frac{1}{T}\int_{0}^{T}a_{j}(x,t,\theta)~d\theta=0.
\end{equation*}
Now to solve \eqref{a.2.10} we expand $ \mathbf{U}_{2} $ in the basis of eigenvectors
\begin{equation}\label{a.2.13}
\mathbf{U}_{2}=\sum_{j=1}^{3}b_{j}(x,t,\theta_{j})\mathbf{R}_{j},\qquad b_{j}=\mathbf{L}_{j}.\mathbf{U}_{2}
\end{equation}
We then use \eqref{a.2.11} and \eqref{a.2.13} in \eqref{a.2.10}, and solve for 
$ \mathbf{U}_{2} $ upon integration we find that,
\begin{equation}\label{a.2.14}
\begin{aligned}
b_{j}(\theta, &\tau)=-\tau[a_{jt}(\theta)+\lambda_{j}a_{jx}(\theta)+k_{j}\mathbf{L}_{j}.
\nabla\mathbf{A}_{0}.\mathbf{R}_{j}\mathbf{R}_{j}a_{j}(\theta)a_{j\theta}(\theta)]\\[5pt]
&-\sum_{m}^{(j)}\sum_{n}^{(j)}k_{m}\mathbf{L}_{j}.
\nabla\mathbf{A}_{0}.\mathbf{R}_{m}\mathbf{R}_{n}
\int_{0}^{\tau}a_{m}[k_{j}^{-1}(k_{m}\theta+\nu_{jm}\zeta)]
a_{n}^{'}[k_{j}^{-1}(k_{n}\theta+\nu_{jm}\zeta)]\,d\zeta+h_{j}(\theta),
\end{aligned}
\end{equation}
where, $\nu_{jm}=k_{j}k_{p}(\lambda_{j}-\lambda_{m}) ,$ superscript $ ' $ represents differentiation with respect to theta, $ h_{j}(\theta) $ is an arbitrary function of integration, and $\displaystyle{\sum_{m}^{(j)}} $ indicates summation avoiding the index $ j .$

In order to obtain equations for the wave amplitudes, the solvability condition  is given by the requirement that the second-order terms $ \mathbf{U}_{2} $ in \eqref{a.2.7} should have sublinear growth in time $ \tau $, that is,
\begin{equation}\label{a.2.15}
\lim_{\tau\to\infty}\frac{1}{\tau}\mathbf{U}_{2}(\xi,\tau)=0.
\end{equation}
This condition ensures that $ \mathbf{U}_{1} $ in \eqref{a.2.7} is a uniformly valid first order approximation for the time of order $ \epsilon^{-1} $, i.e.,
$ \epsilon^2 \mathbf{U}_{2} << \epsilon \mathbf{U}_{1}$ when $ \tau = O(\epsilon^{-1}) .$

In view of sublinearity condition\eqref{a.2.15}, \eqref{a.2.14} reduces to a set of integro-differential equations for $ a_{j} $'s
\begin{equation}\label{a.2.16}
a_{jt}(\theta)+\lambda_{j}a_{jx}(\theta)+
E_{j}a_{j}(\theta)a_{j\theta}(\theta) +
\sum_{m\neq n}^{(j)}\bar{\Gamma}_{jmn}\lim_{T\to\infty}\frac{1}{T}
\int_{0}^{T}a_{m}(\mu_{mnj}\theta+\mu_{mjn}\zeta)
a_{n\zeta}(\zeta)\,d\zeta=0,
\end{equation}
where $\displaystyle{\sum_{m\neq n}^{(j)}} $ denotes the sum over all $ 1 \leq m,n \leq 3\, $ with $ m\neq n,  m \neq j, \,\text{and}\, n \neq j, $
with coefficients given by
\begin{equation}\label{a.2.17}
\setlength{\jot}{5pt}
\begin{aligned}
&E_{j}=k_{j}\mathbf{L}_{j}.\nabla\mathbf{A}_{0}.\mathbf{R}_{j}\mathbf{R}_{j},\\
&\bar{\Gamma}_{jmn}=k_{m}\mathbf{L}_{j}.\nabla\mathbf{A}_{0}.\mathbf{R}_{m}\mathbf{R}_{n},\\
&\mu_{jmn} =  \frac{k_{j}(\lambda_{j}-\lambda_{m})}{k_{n}(\lambda_{n}-\lambda_{m})}.
\end{aligned}
\end{equation}
Upon simplification the interaction terms \eqref{a.2.16} can be written in alternative form
\begin{equation}\label{a.2.18}
a_{jt}(\theta)+\lambda_{j}a_{jx}(\theta)+
E_{j}a_{j}(\theta)a_{j\theta}(\theta) +
\sum_{m < n}^{(j)}\Gamma_{jmn}\lim_{T\to\infty}\frac{1}{T}
\int_{0}^{T}a_{m}(\mu_{mnj}\theta+\mu_{mjn}\zeta)
a_{n\zeta}(\zeta)\,d\zeta=0,
\end{equation}
where $\displaystyle{\sum_{m < n}^{(j)}} $ stands for the sum over all $ 1 \leq m<n \leq 3\, $ with $  m , n \neq j  $ and
\begin{equation}\label{key}
\Gamma_{jmn}=\mu_{jmn} k_{n} \mathbf{L}_{j}.\nabla\mathbf{A}_{0}.\mathbf{R}_{m}\mathbf{R}_{n}
+ \mu_{jnm} k_{m} \mathbf{L}_{j}.\nabla\mathbf{A}_{0}.\mathbf{R}_{n}\mathbf{R}_{m}. \nonumber
\end{equation}
Evaluating the coefficients \eqref{a.2.17}, the wave amplitudes satisfy the following equations
\begin{eqnarray}\label{a.2.19}
\begin{aligned}
&a_{1t}(\theta)-c_{0}a_{1x}(\theta)-
k_{1}\Lambda a_{1}(\theta)a_{1\theta}(\theta) -
k_{2}\Gamma\lim_{T\to\infty}\frac{1}{T}
\int_{0}^{T}a_{2}^{'}   \left(\frac{k_{2}\theta}{2 k_{1}}+\frac{k_{2}\xi}{2 k_{3}}\right)
a_{3}(\xi)~d\xi=0,\\
&a_{2t}(\theta)=0,\\
&a_{3t}(\theta)+c_{0}a_{3x}(\theta)+
k_{3}\Lambda a_{3}(\theta)a_{3\theta}(\theta) +
k_{2}\Gamma\lim_{T\to\infty}\frac{1}{T}
\int_{0}^{T}a_{2}^{'}\left(\frac{k_{2}\theta}{2 k_{3}}+\frac{k_{2}\xi}{2 k_{1}}\right)
a_{1}(\xi)~d\xi=0,\\
\end{aligned}
\end{eqnarray}
where, the coefficients are
\begin{equation}\label{a.2.20}
\begin{aligned}
&\Lambda=c_{0}G ,\qquad c_{0}=\sqrt{\frac{1 + \delta}{1-b}},
\qquad G=\left(1+\frac{\rho_{0}}{c_{0}}c_{\rho 0}\right)=\frac{(\delta + 2)}{2(1-b)},\\
&\Gamma=(c_{\rho_{0}} P_{s_{0}} + c_{0}^{2}c_{s_{0}} + c_{0} P_{s_{0}}/ \rho_{0} - P_{\rho s_{0}}c_{0})=\frac{(1+\delta)^{3/2}}{4(1-b)^{3/2}}.
\end{aligned}
\end{equation}

Equations, \eqref{a.2.19} are two nonlinear Burgers equations for the sound wave amplitudes, coupled by convolution term with kernel given by the derivative of entropy wave amplitude and an entropy equation which implies that the entropy wave is independent of $ \,t\, $ and can be taken as a function of $ (x,\theta) $ determined by the initial conditions.

The quadratically nonlinear self-interaction coefficient of the sound waves $ \Lambda $ is called the parameter of the nonlinearity of the fluid. It is positive for any fluid with normal thermodynamic properties and shows characteristics of classical gas dynamics. For the real gas we considered, it always remains positive for different values of the van der Waal parameter $ b $; whereas the self-interaction coefficient of the entropy wave is zero, which corresponds to the fact that the entropy
wave field is linearly degenerate \cite{MR760229}. The entropy equation has no integral term, i.e., it is not affected by the interaction between the sound waves but it couples the sound waves together. This is explained by the result that the entropy wave is a Riemann invariant of \eqref{a.2.1}.

The nonlinear terms of the integro-differential equations correspond to the Burgers nonlinearity which acts on the acoustic modes and causes steepening of the pressure variation leading to the shock formation, whereas the integral term obtained from the interaction of the acoustic modes with the entropy has a dispersive character \cite{MR1719749} and opposes the shock formation.

Now we focus on the implications of the properties in \eqref{a.2.4} on the solutions and for that matter we write small perturbations of the equilibrium state in the following form
\begin{eqnarray}\label{a11}
\begin{aligned}
&\rho = 1 + \epsilon [ a_{1}(x,t,k_{1}(x+c_{0}t)/\epsilon) + a_{3}(x,t,k_{3}(x-c_{0}t)/\epsilon) + P_{s_{0}} a_{2}(x,t,k_{2} x /\epsilon) ]  + O(\epsilon^2),\\
& v =  c_{0} \epsilon[ a_{1}(x,t,k_{1}(x+c_{0}t)/\epsilon) - a_{3}(x,t,k_{3}(x-c_{0}t)/\epsilon) ]  + O(\epsilon^2),\\
& s =  -c_{0}^{2} \epsilon a_{2}(x,t,k_{2} x/\epsilon) + O(\epsilon^2),\\
\end{aligned}
\end{eqnarray}
it is clear from  \eqref{a.2.4} that the functions $ a_{i}(x,t,\theta) $ are all $ 2\pi$-periodic functions of their arguments with phases phases $ \theta_{1}, \theta_{2}, \text{and}~ \theta_{3} $.

Next, using the Fourier series expansions and periodicity of $ a_{i}$'s it is apparent that the resonant integral terms in \eqref{a.2.19}(a)  vanish unless  $ \displaystyle {\frac{k_{2}}{2 k_{3}}} $   is a rational number. Similarly, limiting integral terms in \eqref{a.2.19}(c) vanish unless  $ \displaystyle {\frac{k_{2}}{2 k_{3}}} $ is a rational number.
Alternatively, the resonance condition for  periodic solutions is given by \cite{MR1266838}, 
\begin{eqnarray}\label{a.2.21}
\begin{aligned}
&k_{2} = r k_{1} + s k_{3},\quad \omega_{2} = r\omega_{1} + s\omega_{3},\quad \\
&\omega_{1} = -c_{0}k_{1}, ~~\omega_{2} = 0,~~\omega_{3} = c_{0}k_{3},
\end{aligned}
\end{eqnarray}
where $ r $ and $ s $ are rational numbers.
Equations \eqref{a.2.21} imply
$ \displaystyle{\frac{k_{2}}{2k_{1}}=r,\,\,\frac{k_{2}}{2k_{3}}=s.} $

We consider the strongest interaction between the the waves, which occurs when the fundamental harmonics of both wave resonate directly. In our case, this corresponds to the situation $ m=1, ~ n=1 ,$ so that there is a direct resonance between the  fundamental harmonics of the entropy and sound waves.

In view of the above discussions and assuming that there are no spatial modulations (the wave amplitudes do not depend on $ x $), the system \eqref{a.2.19}  reduces to the following form
\begin{eqnarray}\label{a.2.50}
\begin{aligned}
&a_{1t}(\theta)-
\Lambda a_{1}(\theta)a_{1\theta}(\theta) -
\frac{\Gamma}{2\pi}
\int_{0}^{T}a_{2}^{'} \left(\theta+\xi\right)
a_{3}(\xi)~d\xi=0,\\
&a_{3t}(\theta)+
\Lambda a_{3}(\theta)a_{3\theta}(\theta) +
\Gamma\frac{1}{2\pi}
\int_{0}^{T}a_{2}^{'} \left(\theta+\xi\right)
a_{1}(\xi)~d\xi=0.\\
\end{aligned}
\end{eqnarray}
Untill now we have applied only the periodic properties in \eqref{a.2.4} for \eqref{a11} but not the symmetries ($ \rho $ is even etc,). It is apparent that symmetric properties are equivalent to the following identity between $ a_{1} $ and $ a_{3} $:
$ a_{1}(t, \theta) $ = $ a_{3}(t,-\theta) $ and  $ \sigma_{2} $ is an even function of phase $ \theta.$

We now introduce the common function 
\begin{equation}\label{a.51}
\sigma (t,x)= a_{1}(t,-x)=a_{3}(t,x),
\end{equation}
and  the function 
\begin{equation}\label{a.52}
K(x)= a_{2}'(t,x)/2\pi,
\end{equation}
where we have switched the notation for phase arguments, using x instead of $ \theta $; this should not lead to any confusion as there is no dependence of amplitudes on the original space variable.
$ K(x) $ is an odd $2\pi$-periodic function of its argument, which can be determined from the initial conditions. Hence, it has a Fourier sine series expansion 
\begin{equation}\label{a.53}
K(x)=\sum_{n=1}^{\infty} A_{n} \sin (nx).
\end{equation}
Finally, equations \eqref{a.51} and \eqref{a.52} imply that the whole system \eqref{a.2.50} collapses into the single equation
\begin{equation}\label{a.2.24}
\begin{aligned}
&\frac{\partial}{\partial t} \sigma(x,t)+
\Lambda \frac{\partial}{\partial x}\left[\frac{1}{2}\sigma^2(x,t)\right] + \Gamma
\int_{0}^{2\pi}K(x-y)
\sigma(y,t)~dy=0.\\
\end{aligned}
\end{equation}
This equation is the main object of our study in this chapter. The solution $ \sigma(x,t) $
is $2\pi$-periodic in phase x and has zero mean (for all time t). We observe that the equation is compatible  with these requirements because $ K $ has zero mean. In this chapter we limit our discussion to the simplest possible kernel $ K(x)=\sin x $ that represent a periodic wave.



\subsection{Traveling wave solutions for the system}
We have seen in the last section  that in the absence of space modulations (the wave amplitudes are independent of $ x $) our problem of gas motion inside a duct is reduced to the system \eqref{a.2.50}.
In the ideal gas background, a system similar to this was studied extensively by Majda, Rosales, Schnobek \cite{MR975485} . They found the analytical and numerical solution for both local and nonlocal systems obtained by choosing the appropriate kernels  $ K $ in \eqref{a.2.50}, and found out an interesting property that they possess a family of traveling wave solution. Pego \cite{MR975486} and Gabov \cite{MR510246} calculated an exact analytic expression for these solutions. For the real gas case, with sinusoidal entropy distributions, the explicit formulation of the traveling wave solutions is given by
	Suppose that $  K(x) = \sin x $. Then two families of traveling waves of \eqref{a.2.50}
	in the form of \eqref{a.2.26} with zero mean, exist. Fixing a choice of $ \delta =\mp 1$ for any $ \alpha \in [0,1],  (u_{1},u_{2} ) $ and traveling speed $ s $ may be given by
	\begin{equation}\label{a.2.26}
	\setlength{\jot}{10pt}
	\begin{aligned}
	&u_{1}(\psi)=
	\delta \left(\frac{\Gamma}{\Lambda} \right)\gamma \sqrt{1+\alpha \cos\psi} + \frac{s}{\Lambda}, \\
	&u_{2}(\psi)=
	\delta \left(\frac{\Gamma}{\Lambda} \right)\gamma \sqrt{1+\delta\, \alpha \cos\psi} + \frac{s}{\Lambda}, \\
	&s=-\delta \left(\frac{\Gamma\gamma}{2\pi} \right)\int_{0}^{2\pi} \sqrt{1+\alpha \cos y}~dy, \\
	\end{aligned}
	\end{equation}
	where $$ \gamma = \gamma(\alpha) =\frac{2}{\alpha}\int_{0}^{2\pi}\cos y \sqrt{1+\alpha \cos y}~dy, \quad\text{for}\quad \alpha\neq 0 .$$
	These traveling waves exist up to a maximum amplitude corresponding to $ \alpha=1 $ and the maximum amplitude waves are cusped.
\begin{proof}
	
	We look for traveling waves of the following form
	\begin{equation}\label{a.2.27}
	\sigma_{1}=u_{1}(x-st),\qquad\sigma_{2}=u_{2}(x-st).
	\end{equation}
	Using\eqref{a.2.27} in \eqref{a.2.26} with $ K(x) = \sin x $ and upon integration, we get
	\begin{equation}\label{a.2.28}
	\setlength{\jot}{10pt}
	\begin{aligned}
	&-su_{1}+
	\frac{1}{2}\Lambda u_{1}^2 = \Gamma
	\int_{0}^{2\pi}\cos(z-\zeta)
	u_{2}(\zeta)~d\zeta+k_{1},\\
	&-su_{2}+
	\frac{1}{2}\Lambda u_{2}^2 = \Gamma
	\int_{0}^{2\pi}\cos(z-\zeta)
	u_{1}(\zeta)~d\zeta+k_{2},\\
	\end{aligned}
	\end{equation}
	where $ k_{1}, k_{2} $ are constants of integration. Substituting $ w_{i} = u_{i} - s/\Lambda$ in \eqref{a.2.28} the resulting system for some non-negative 
	$ \gamma_{1},\gamma_{2} $ reduces to the following system on using cosine expansion
	\begin{equation}\label{a.2.30}
	\setlength{\jot}{10pt}
	\begin{aligned}
	&w_{1}^2 = \left(\frac{4\pi\Gamma}{\Lambda}\right)
	[C_{w2}\cos(z)+S_{w2}\sin(z)] + \left(\frac{\Gamma\gamma_{1}}{\Lambda}\right)^2,\\
	&w_{2}^2 = \left(\frac{4\pi\Gamma}{\Lambda}\right)
	[C_{w1}\cos(z)+S_{w1}\sin(z)] + \left(\frac{\Gamma\gamma_{2}}{\Lambda}\right)^2,\\
	\end{aligned}
	\end{equation}
	where, $  C_{w1},C_{w2},S_{w1} ~\text{and}~ S_{w2}  $ are the fourier coefficients of $w_{1} ~\text{and}~ w_{2}, $
	\begin{equation}\label{a.2.31}
	\setlength{\jot}{10pt}
	\begin{aligned}
	&C_{w1} = \frac{1}{2\pi}
	\int_{0}^{2\pi} \cos(z)w_{1}(z)\,dz,\qquad 
	C_{w2} = \frac{1}{2\pi}
	\int_{0}^{2\pi} \cos(z)w_{2}(z)\,dz,\\
	&S_{w1} = \frac{1}{2\pi}
	\int_{0}^{2\pi} \sin(z)w_{1}(z)\,dz,\qquad 
	S_{w2} = \frac{1}{2\pi}
	\int_{0}^{2\pi} \sin(z)w_{2}(z)\,dz,\\
	\end{aligned}
	\end{equation}
	Since \eqref{a.2.28} is translation invariant, replacing $  (w_{1},w_{2})(x) $ by an appropriate translate such that $ S_{w2}=0,$ which implies $ w_{1} $ is an even function and it follows that $S_{w1}=0.$  
	
	Using \eqref{a.2.30} we can find that for some choice of 
	$ \,\delta_{i},\,\delta \in {\{-1,1\}}, $
	\begin{equation}\label{a.2.32}
	w_{1}=
	\delta_{1} \left(\frac{\Gamma}{\Lambda} \right)\gamma_{1} \sqrt{1+\alpha_{2} \cos z}, \qquad
	w_{2}=
	\delta_{2} \left(\frac{\Gamma}{\Lambda} \right)\gamma_{2} \sqrt{1+\delta\, \alpha_{1} \cos z}, 
	\end{equation}
	where  \begin{equation*}\label{key}
	\alpha_{1}=\left(\frac{4\pi\Lambda}{\Gamma} \right)\gamma_{2}^{-2}C_{w1}, \qquad
	\alpha_{2}=\left(\frac{4\pi\Lambda}{\Gamma} \right)\gamma_{1}^{-2}C_{w2}, 
	\end{equation*}
	\noindent
	using \eqref{a.2.32} in \eqref{a.2.30} we obtain 
	\begin{equation}\label{a.2.33}
	\setlength{\jot}{10pt}
	\begin{aligned}
	&\alpha_{1}\gamma_{2}^{2} = 2\delta_{1}\gamma_{1}\int_{0}^{2\pi}\cos z \sqrt{1+\alpha \cos z}~dz,\\
	&\alpha_{2}\gamma_{1}^{2} = 2\delta_{2}\gamma_{2}\int_{0}^{2\pi}\cos z \sqrt{1+\delta \,\alpha \cos z}~dz.\\
	\end{aligned}
	\end{equation}
	Now, consider the function
	\begin{equation}\label{a.2.34}
	P({\alpha})=\int_{0}^{2\pi}\cos x \sqrt{1+\alpha \cos x}~dx=
	\int_{-1}^{1}\sqrt{1+az}\,\frac{z\,dz}{1-z^{2}}.
	\end{equation}
	It is easy to see that $ P $ is an odd function of $  \alpha $ and positive for $ \alpha > 0, $ then from \eqref{a.2.33} it follows that $ \delta=\delta_{1}=\delta_{2}.$ 
	
	Now it only remains to show  that $ \alpha_{1} = \alpha_{2}, $ which implies $ \gamma_{1} = \gamma_{2} $, then using the property that $ u_{1},u_{2}$ have zero mean we get the required expression \eqref{a.2.26}. 
	To meet this objective we define a function 
	\begin{equation}\label{a.2.35}
	Q({\alpha})=\int_{0}^{2\pi} \sqrt{1+\alpha \cos x}~dx=
	\int_{-1}^{1}\sqrt{1+\alpha z}\,\frac{dz}{1-z^{2}}.
	\end{equation}
	which implies
	\begin{equation}\label{a.2.36}
	\gamma_{1}Q(\alpha_{2}) = \gamma_{2}Q(\delta\alpha_{1}).
	\end{equation}
	As, $ Q (\alpha) $ is even and decreasing for $ \alpha $ positive, using \eqref{a.2.33} and \eqref{a.2.35} we get
	\begin{equation}\label{key}
	\frac{\alpha_{1}P(\alpha_{1})}{Q(\alpha_{1})^{3}}=
	\frac{\alpha_{2}P(\alpha_{2})}{Q(\alpha_{2})^{3}}
	\end{equation}
	except for the trivial case $ u_{1} = u_{2} = 0. $
	Since the function $ \displaystyle{\frac{\alpha P(\alpha)}{Q(\alpha)^{3}}}  $ is increasing for $  \alpha \geq 0, $ it follows that $  \alpha_{1} = \alpha_{2}. $ which completes the proof since all the results of the proposition follow from 
	$ \gamma=P(\alpha)/\alpha$.
\end{proof}

\section{Traveling wave solutions}\label{trav}
In this section, we study the traveling wave solution for the single evolution equation \eqref{a.2.24} and draw attention to some of their important characteristics in the view of real gas settings. Following Pego \cite{MR975486} we obtained a one parameter family of traveling wave solutions for the system \eqref{a.2.50}, which (except for phase shift) in the case of our single transport equation (for a single mode kernel $ K(x)=\sin x $ ) reduces to the following form
\begin{equation}\label{a.2.37}
\sigma (x-st) = \frac{s}{\Lambda} + \left(\frac{\Gamma}{\Lambda} \right)\gamma \sqrt{1+\alpha \cos(x-st+\phi)} ,
\end{equation}
where $\, \phi \,$ is a phase shift; $ s=s(\alpha) $ is the speed of wave propagation, obtained by assuming that $ \sigma $ has zero mean and is given by
\begin{equation}\label{a.2.38}
s(\alpha)=-\left(\frac{\Gamma\gamma}{2 \pi} \right)\int_{0}^{2\pi} \sqrt{1+\alpha \cos y}~dy,
\end{equation}
and $ \gamma(\alpha) $ has an expression
\begin{equation}\label{a.2.39}
\gamma(\alpha)=\frac{2}{\alpha} \int_{0}^{2\pi}\cos y \sqrt{1+\alpha \cos y}~dy.
\end{equation}
Evidently, the parameter $ \alpha $ is confined to the range $ 0\leq|\alpha|\leq 1 $
and a change of sign in $ \alpha $ is equivalent to a phase shift of $ \pi  $ in the wave also for all values of $ \alpha,~ s<0 ~\text{and}~ \gamma(\alpha)>0.$

The amplitude of the wave (the half of the distance between the maximum and the minimum) depends on the parameter $ \alpha $ and $ \gamma(\alpha) $ and is given by
\begin{equation}\label{a.2.40}
A=\left(\frac{\Gamma}{2\Lambda} \right)\gamma (\sqrt{1+|\alpha|}-\sqrt{1-|\alpha|}).
\end{equation} 

\begin{figure}[!t]
	\centering
	\includegraphics[width=1.02\textwidth]{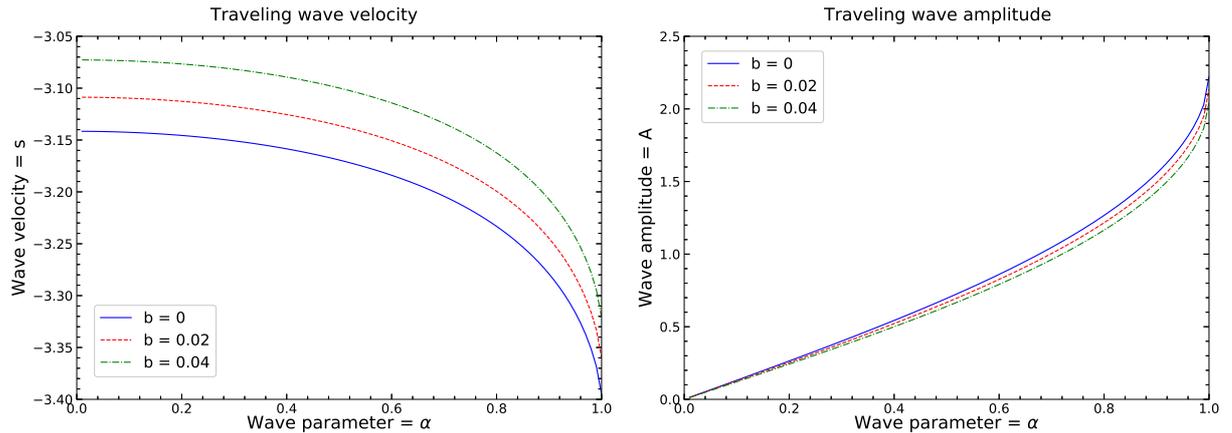}
	\caption{Wave velocity (upper left) and wave amplitude (upper right) variation for various values for $b=0,\,b=0.02,\,~\text{and}~ b=0.04 $ as a functions of parameter $ \alpha $ for $ 0 \leq \alpha \leq 1 .$}
	\label{100}
\end{figure}

Three sets of computations are performed for each of the traveling wave parameters above for different values of the van der Waal parameter $ b $ namely; $ b=0,\, b=0.02,\, b=0.04 $
with fixed value of the $ \delta=0.4 $ for air; here $ b=0 $ corresponds to the ideal gas behavior. The effect of van der Waal parameter $ b $ on traveling waves solutions of equation \eqref{a.2.24} is shown in Figure (\ref{100})

\begin{figure}[!t]
	\label{101}
	\centering
	\includegraphics[width=1.0\textwidth]{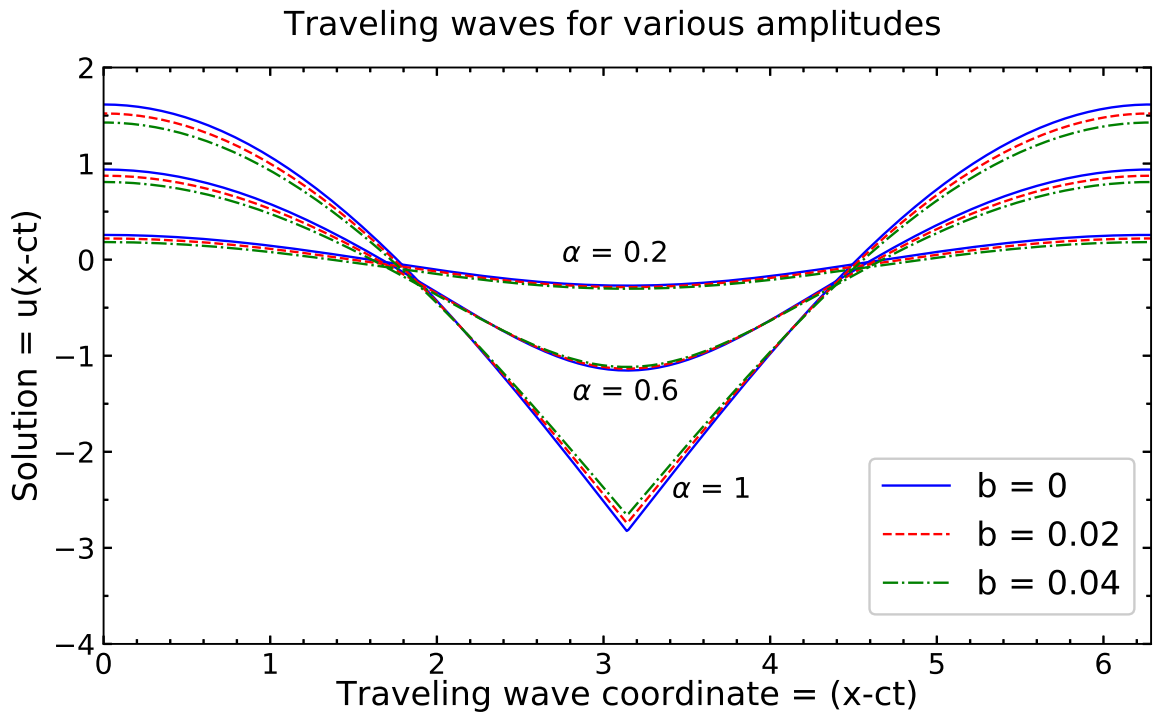}
	\caption{Three families of traveling wave profiles for $b=0,\,b=0.02,\,~\text{and}~ b=0.04 $ for three values of  parameter $\alpha=0,\,0.6,\,~\text{and}~ 1$ with fixed  $\delta=0.4$. The first two waves $( \alpha = 0.2, 0.6) $ of each family are smooth, whreas the third one $( \alpha=1) $ exhibits a maximum amplitude with corner singularity.}
	\label{fluxb0}
\end{figure}


In Fig.(\ref{100})(a), the variation of wave velocity $ s $ with parameter $ \alpha $ is depicted for three values of real gas parameter. It is seen that wave velocity decreases with  increase in the parameter $ \alpha $ in all the three cases, whereas an increase in van der Waal parameter causes the overall increase in the value of the velocity $ s $. 

In Fig.(\ref{100})(b), the change in the amplitude $ A $ with parameter $ \alpha $ is shown for different values of the real gas parameter. It is noticed that the wave amplitude increases with  increase in the  parameter $ \alpha $ while it decreases with an increase in the van der Waal parameter $ b $. 


Finally, in Figure (\ref{fluxb0}) for each value of the van der Waal parameter three typical wave profiles are shown corresponding to the values of parameter $ \alpha = 0.2,\,0.6,\,\text{and}\,1.$
In each case, for a fixed value of real gas parameter the traveling wave solution forms a family which ranges from trivial solution $ \sigma\equiv 0, $ at $ \alpha=0 $ to a limiting maximum amplitude for $ \alpha =\mp 1, $ of which three members are shown.
The traveling waves are smooth (in fact analytic) function of the argument 
$ \psi = x-st+\phi $ for $ |\alpha|< 1$ but as soon as $ \alpha $ reaches the critic value
$ \alpha= \mp1 $ the solution display non smooth structure and develops a corner singularity in all the three cases. In fact our numerical result shows that presence of corners is distinctive for the solution of \eqref{a.2.24} in ideal gas as well as in real gas cases.  

	The traveling wave solutions are periodic in space as well as in time with period 
	$ 2\pi/|s| $. They are limited to a range for $ |\alpha|<1 $ and our numerical experiment indicated that they are the only (stable) time periodic solutions of equations.

\section{Numerical methods and solutions}\label{num}

This section is devoted to the study of the numerical solutions of the evolution equation \eqref{a.2.24}  obtained by performing a long-time numerical integration. In the ideal gas background (b=0), for the transport equations like \eqref{a.2.24}, the existence of NBAT (Non Breaking For All Times) solutions in the form of traveling waves have been shown both analytically and numerically \cite{MR975485}, moreover, the traveling wave for \eqref{a.2.24}  corresponds to the nonlinear standing acoustic waves of the Euler equations \eqref{a.2.1}.  In addition to the traveling waves, Rosales and Shefter (\cite{MR1719749},\cite{MR2716775}) had shown numerically the existence of another family of NBAT solutions that are quasiperiodic in time and are globally attracting in nature for the long-time evolution of the equations.
These new solutions unite both behaviors of having sharp spikes when the corner singularities emerge and being smooth for the rest of the time. Whereas, traveling wave solutions are either smooth or have a corner discontinuity for all time.

We investigated the existence of the NBAT solutions in the light of real gas scenario for equation \eqref{a.2.24}  and discussed the qualitative and quantitative changes in the behavior of NBAT solution for equation \eqref{a.2.24}  with van der Waal parameter $ b. $

First, we give a description of the numerical scheme used in the experiment. As mentioned, in the last section, the presence of the corners in the solution of the equation \eqref{a.2.24}  is a distinctive feature, as noted in the solution  \eqref{a.2.37} for the maximum amplitude traveling wave for  $ \alpha=\mp1, $ and they remain nonsmooth for all time. Besides, shocks are also present in the solution (at least for a short period of time) after
the initial data starts evolving. In our experiment, even for smooth initial data, we almost always find both shocks and weak discontinuities during the evolution of the equation \eqref{a.2.24}.

In our case, $ \Gamma, \Lambda > 0, $ the procedure of formation of strong shocks is similar to that for the inviscid Burger equation 
$$ \displaystyle{\sigma _{t}(x,t)+\left(\frac{\sigma^2(x,t)}{2}\right)_{x}=0,} $$ 
nonlinear term $\, \sigma\, \sigma_{x}\, $ produces the steepening in the equation \eqref{a.2.24} that leads to breaking of wave but the equation has linear convolution term $ K*\sigma $ 
with a weakly dispersive character. It opposes the wave breaking developed by the nonlinear term. The existence of the never breaking solution with corner singularity is the result of the balance between these two forces acting against each other.

Finite difference and spectral methods are efficient on smooth solutions however, they break down in presence of discontinuities as described above. The adequate representation of the discontinuities is very important. If not resolved properly it produces numerical noise which diffuses into the smooth regions and eventually contaminates the solution. Taking into account the long-time integration involved in the experiment it becomes impossible to tell whether the observed solution is real or just accumulated numerical errors. Hence the scheme must satisfy the requirement of safely handled singularity of both strong (shock) and weak (corner) type.

In hyperbolic equations such as \eqref{a.2.24}, directions of propagation of information are governed by the characteristics. The numerical scheme must take into account these directions otherwise information from the physically immaterial domain may be used by the scheme and produces spurious oscillation and other instabilities near to the discontinuities.

In our experiment we perform three sets $ (b=0,0.02,0.04) $ of numerical integration of the equation \eqref{a.2.24} for long time ($ 100-1000 $ acoustical periods) with arbitrary shape initial data. Hence, the numerical scheme is required to capture high gradient regions without spurious oscillations as well as resolve smooth features of the solution with high-order accuracy.
In order to take account of all these requirements, we employ a fractional step approach \cite{MR2731357}, in which we split the equation \eqref{a.2.24} into two subproblems that can be solved independently. First, we split the spatial operator,
\begin{equation}\label{2.17}
\frac{\partial \sigma (x,t)}{\partial t} = P [\sigma(x,t)] + Q [\sigma (x,t)],
\end{equation}  
where $$ P [\sigma(x,t)] = - \Lambda \frac{\partial}{\partial x} (\sigma^{2}(x,t)), \quad 
Q[\sigma(x,t)] = - \Gamma ( K*\sigma) .$$
On discretization of the time variable into equal time steps $ [t_{n-1},t_{n}], $ we obtain the following semi-continuous equation an approximation of \eqref{2.17}
$$\frac{\sigma^{n+1}(x,t)-\sigma^{n}(x,t)}{\Delta t } = P [\sigma(x,t)] +  Q [\sigma (x,t)],  $$
where, $ \sigma^{n} = \sigma (x,n\Delta t) $ is the  $ n$-th time step numerical approximation for a solution and $ \sigma^{0} = \sigma (x,0) $ is the initial condition.
Now, the idea of the fractional step is to combine the simpler problems in an alternating manner in order to approximate the solution of the complete problem. In order to yield higher order of accuracy we use Strang splitting \cite{MR0235754}, which uses a slight modification in the split idea. We solve first sub-problem  over half a time step,
$$ \sigma^{'n}(x,t) = \sigma^{n}(x,t) + \bar{Q} [\sigma^{n}(x,t)]\frac{\Delta t}{2}, $$
then, we use the result data for a full time step in another subproblem 
$$ \sigma^{*n}(x,t) = \sigma^{'n}(x,t) + \bar{P} [\sigma^{'n}(x,t)]\Delta t, $$
and finally take another half step on first subproblem
$$\sigma^{n+1}(x,t) = \sigma^{*n}(x,t) + \bar{Q} [\sigma^{*n}(x,t)]\frac{\Delta t}{2}, $$
where, $ \bar{P} $ and $ \bar{Q} $ are discrete operators corresponding to $ P $ and $ Q $.
Since each intermediate step involve only one spatial operator, we can use different discretization scheme for incorporating different behavior of the operators $ P $ and 
$ Q $. We aim to achieve for the overall second-order accurate scheme, for which it is sufficient to construct second-order approximation in both $ \bar{P} $ and $ \bar{Q}$ as shown by Strang  \cite{MR1275838} .

In the nonlinear subproblem of the equation \eqref{2.17} we encounter discontinuities in the solutions hence we use schemes that are specifically built to treat the discontinuous solutions. We used a second-order Godunov scheme \cite{MR0119433, MR1925043} for the problem. It performs excellently on the strong discontinuities and does not produce spurious oscillations. To achieve the second-order accuracy and good resolution of corners we use the UNO reconstruction algorithm \cite{COLELLA1984174, MR1419777} in the Godunov scheme. To resolve shocks effectively we use an adjustable grid approach \cite{MR707200, MR881365} which allows grid points to move along as the shock passes through them.

Next, we consider the linear convolution part of the splitting scheme. To preserve second-order accuracy of the overall scheme we used the trapezoidal rule to perform the numerical integration and the temporal variable is discretized using an explicit second-order Runge-Kutta scheme \cite{MR0368379}.

\subsection{Evolutionary stages for general initial conditions } 
To study the long-time behavior of equation \eqref{a.2.24} we conducted substantial numerical experiments in the real gas setting with arbitrary $ 2\pi $ periodic function of various shapes as initial data.  The influence of real gas background on long term behavior on these solutions has been studied.
These computations gave the confirmation of the existence of the NBAT solution in the real gas environment as well. These solutions form a family of globally attracting general solution that arises from arbitrary initial data. 

We carried out three-set of numerical experiments for three different values of the real gas parameter $ b $ viz.  $ b=0,\, b=0.02,\,\text{and}\, b=0.04 $ with \eqref{a.2.24} (with the kernel 
$ K(x)=\sin x $) and the result demonstrates that the evolution of the solution is representative of almost all initial conditions. Broadly it is composed of the following two stages

\subsubsection{1. Initial temporary stage.} 
\begin{figure}[!t]
	\centering
	\includegraphics[width=1.00\textwidth]{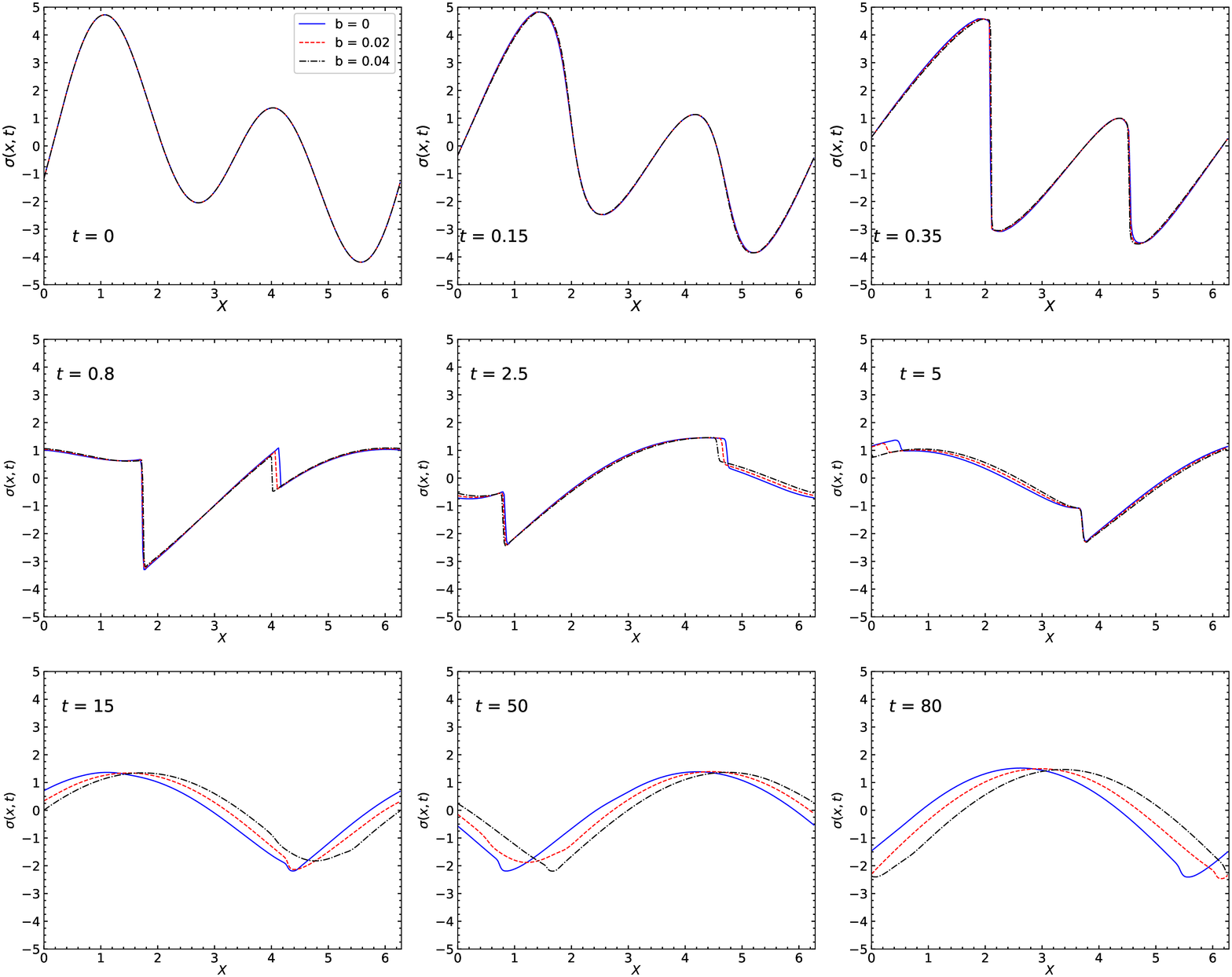}
	\caption{ Evolution of three families of wave profiles corresponding to  $b=0,\,b=0.02,\,~\text{and}~ b=0.04 $ with $\sigma(x, 0)=2 \sin x +3 \cos (2x-2) $ with fixed  $\delta=0.4$. The horizontal and vertical axis are $ 0<x<2\pi, -5 < \sigma(x,t) < 5 $, respectively. From top to bottom and left to right  : initial conditions and $ \sigma (x, t) $ at time $ t =\,0.15,\,0.35,\,0.8,\,2.5,\,5,\,15,\,50,~\text{and}~80\,.$ Shocks first appears around $ 0<t<0.35 .$ At first, they build up, as the nonlinearity steepens the solution, but by $ t=0.8$ they are decaying in strength and by $ t=5 $ they are weak. The shocks are disappeared by $ t=50 $ and no new shock is formed after; however, notable acoustic wave amplitude is present in the solution.}
	\label{fig131}
\end{figure}

\begin{figure}[!t]
	\centering
	\includegraphics[width=1.00\textwidth]{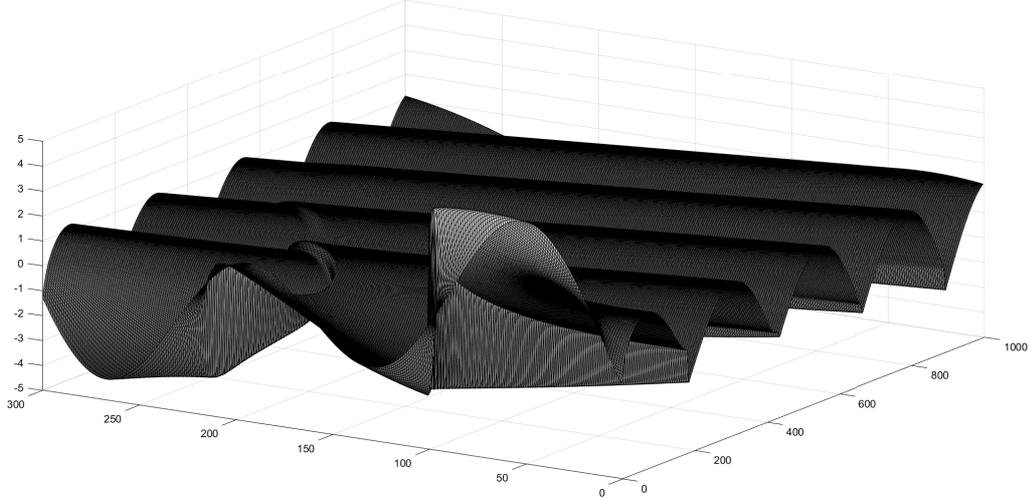}
	\caption{Front view of the evolution of the solution in a physical space-time  for $ b=0.02 $ with $\sigma(x, 0)=2 \sin x +3 \cos (2x-2) $ and fixed  $\delta=0.4$ from $ 0<t<10. $ Initially, shocks have formed but after some time they disappear even though solution has a substantial wave amplitude.}
	\label{fig4}
\end{figure}

\begin{figure}[!t]
	\centering
	\includegraphics[width=1.05\textwidth]{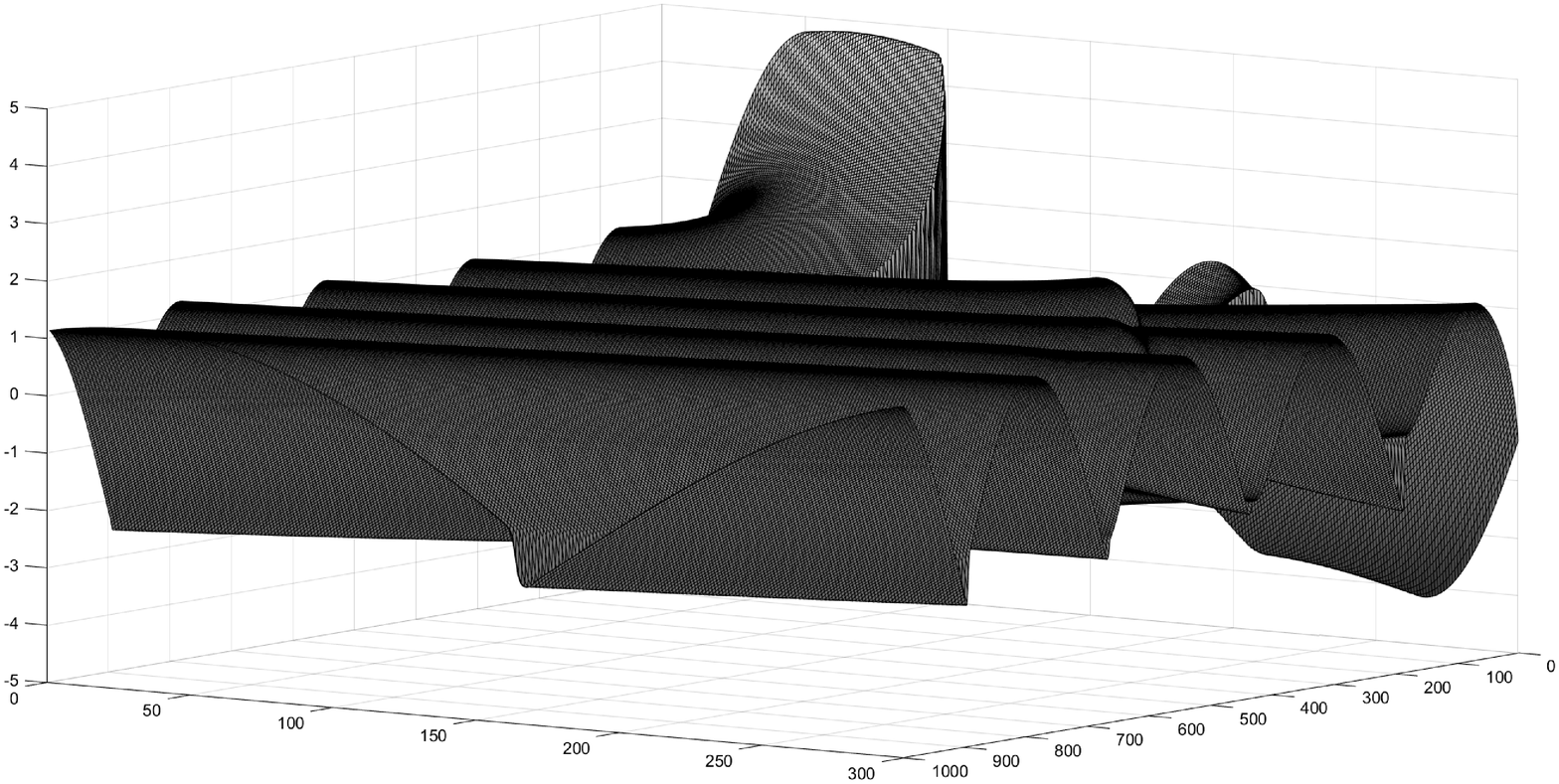}
	\caption{Rear end view of the evolution of the solution in  physical space-time for $ b=0.02 $ with $\sigma(x, 0)=2 \sin x +3 \cos (2x-2) $ ,fixed  $\delta=0.4$, $ 0 <x < 2\pi$, and $ 0<t<10. $ The shocks have disappeared but wave amplitude is nonzero in the solution.}
	\label{fig5}
\end{figure}

In all the three cases just after the evolution begins steepness of the wave profile increases and shocks appear similar to the case of nonviscous Burgers equation. Energy dissipation occurs due to the presence of shocks, and they become weaker with time. After some duration of time (generally less than a hundred period) the shock becomes so weak that it is virtually impossible to detect them by the current numerical scheme, hence they can be effectively ignored. However, the energy is still nonzero at this stage. There are no observable differences in the evolution of the solution for different values of the real gas parameter so far. For continuous initial data, a typical example of this type of evolution is shown in Figures (\ref{fig131}),(\ref{fig4}),(\ref{fig5}).

\subsubsection{2. Final limiting stage.}
After the shocks vanish from the solution profile, it attains a limiting stage. In all the experiments we performed, no more shocks are formed during the final stage. However, the solution still has significant pressure variations. We call these solutions attained by this transformation mechanism the limiting solution and examine their nature.  

Shafter and Rosales \cite{MR1719749, MR2716775}  studied these solutions in the ideal gas background. In our experiment, we found that these solutions are also present in the real gas case. The behavior of these solutions are quite different from the nonlinear counterpart of equation \eqref{a.2.24}, i.e., an equation similar to the Burgers equation, where all the energy is dissipated by the shocks, or completely dispersive equation such as KdV equations where shocks do not form at all. It is due to the weakly dispersive character introduced by the convolution term in the equation. Weakly dispersive means dispersiveness becomes less important  as the wavelength goes to zero, this is the reason dispersive term can balance the nonlinear term and prevent the shock formation only for a restricted range of amplitudes. It can be seen 
by considering the linearization of the equation \eqref{a.2.24} for the small values of $ \sigma(x,t) $ of the form
\begin{equation}\label{a.2.74}
\begin{aligned}
&\frac{\partial}{\partial t} \sigma(x,t)+
\Gamma
\int_{0}^{2\pi}K(x-y)
\sigma(y,t)~dy=0.\\
\end{aligned}
\end{equation}
This is a translation invariant linear equation, hence it can be solved by separation of variables and has elementary solutions of the form $ \sigma(x,t) = A e^{i(kx-\omega t)}$. Since we are looking for $ 2\pi $-periodic solutions with zero mean, it restricts the wave number to $ k = \pm 1,\pm 2,\pm 3,\ldots $ It is easy to see that for this being a solution, $ \omega  $ satisfies the dispersion relation 
\begin{equation}\label{key}
\omega = - \Gamma i \int_{0}^{2 \pi} K(\xi) e^{-ik\xi} d\xi = - \Gamma 2 \pi i \widehat{K}(k)
\end{equation}
where $ \widehat{K}(k) $ is the $ k $th Fourier coefficient of $ K $. Since $ K $ is an odd function, $ \omega $ is real. Moreover, $ \omega $ is not a linear in $ k $, because $ K $ is a function. Therefore, the equation is dispersive.
It is noticeable  that the wave number $ \omega $ vanishes as $ |k|\to \infty $. Thus, for short
waves the convolution term does not have a strong effect and  the strength of this dispersion to smooth out the solution does not increase as the waves become steeper, this is what we mean by weak dispersion.  It will not always balance the nonlinear term, as would be the case in a strongly dispersive (Kdv equation) equation. When the steepening 
effects are not too strong, a weak dispersion can be effective in  keeping the nonlinear advective term balanced. However, generally it cannot prevent the formation of shocks.

	If the solution is not in resonance with the kernel of the integral term we do not get the NBAT solutions. Since in that case, the contribution from the integral term is zero thus; the equation \eqref{a.2.24} reduces to the inviscid burger equation. Then as typical of the Burgers equation, the energy decreases asymptotically to zero.

\subsubsection{Attracting solution structure} 
The limiting solutions obtained as a result of the numerical experiment with different initial condition show many common characteristics and can be treated as a family. When the initial conditions have sufficiently large energy  \cite{MR1719749, MR2716775} , the limiting solutions obtained are maximum amplitude traveling waves  $ (\alpha=1 )$ and are given by
\begin{equation}\label{key}
\sigma = \frac{s(1)}{\Lambda} + \left(\frac{\Gamma}{\Lambda} \right)\gamma(1) \sqrt{1+\alpha \cos(x-s(1)t+\phi)} ,
\end{equation}
for some phase shift $ \phi,\,$ where $ s(1)=-32\Gamma/3\pi \,,$ and $\, \gamma(1)=8\sqrt{2}/3. $
\begin{figure}[!t]
	\centering
	\includegraphics[width=1.00\textwidth]{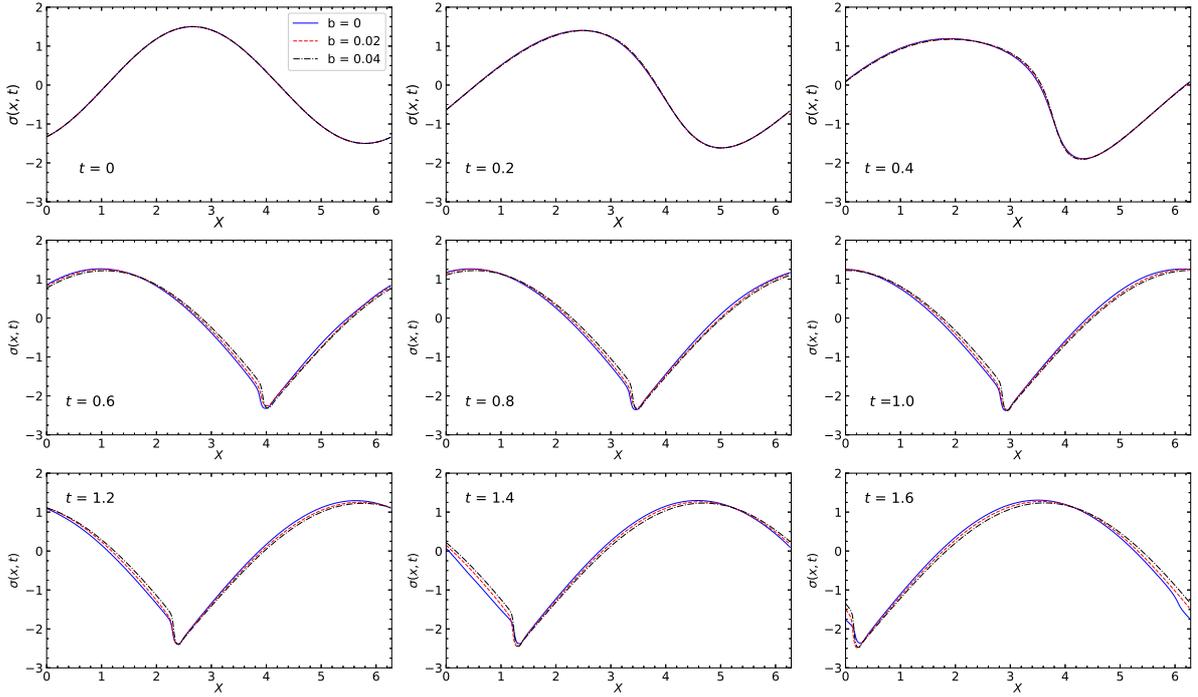}
	\caption{Evolution of three families of a typical attracting solution following from the initial condition $ \sigma(x,0)=1.5\,\sin x $ illustrating the formation of a corner discontinuity near the trough of the wave formed at $ t=0.6 $ and moves across the bottom to the other side. The horizontal and vertical scale are $ 0\leq x \leq 2 \pi\, \text{and}\, -3 \leq \sigma(x,t) \leq 2$, respectively. Time proceeds from top to bottom and left to right in steps $ 0.2 $ in each calculation $ \delta=0.4 $ is fixed.}
	\label{fig 6}
\end{figure}
\begin{figure}[!t]
	\centering
	\includegraphics[width=1.00\textwidth]{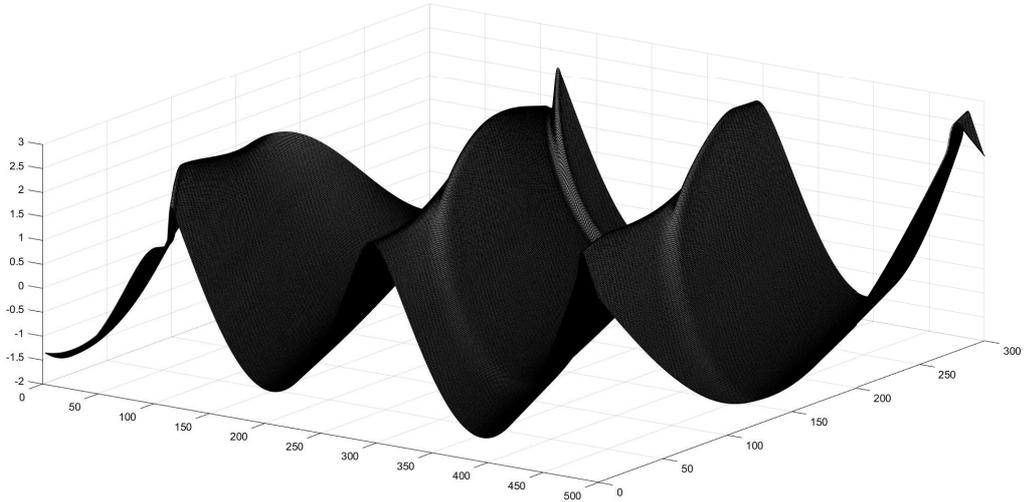}
	\caption{Space-time (top view) of an attracting  solution. The tracks formed by the movement of the corners along the bottom and causes the periodic shape change. We have depicted the negative of the solution to make it more evident. }
	\label{fig7}
\end{figure}

When the initial condition does not have sufficiently large energy, experiments suggest the existence of a completely different class of solutions. These solutions are periodic like traveling wave solution but unlike traveling waves, they change their shape because an oscillatory perturbation is superimposed on them. Their typical structure is shown in Figs. (\ref{fig 6}),(\ref{fig7}) in ideal as well as real gases.
In all the cases the overall structure of the limiting solutions can be explained in a manner similar to that of a smooth traveling wave with an oscillation per wavelength, but there is a periodic shape change with the movement as explained next.
In all the wave profiles, with the start of each cycle, there is a formation of a single corner on the front part near to the trough then corner moves to the rear side along  the wave profile and vanishes. After the formation of the corner, the solution remains smooth for some period then the cycle starts again.

It is noticeable that in spite of not being smooth for all time the solution never breaks and forms no shocks. This behavior is quite different from that of traveling waves which either have slope discontinuity or are smooth for all time as described in section (\ref{trav}). These solutions show a combination of the two behavior showing spikes during the presence of corner singularities or otherwise being smooth for all time.
A repetitive pattern is clearly seen in the behavior of these solutions in the experiment. Numerical results suggest the shape oscillations are periodic in time. Since the period of translation of shape is combined with the period of shape oscillation these solutions are quasiperiodic in time having two periods. 

	The quasiperiodic solution corresponds to nonlinear almost standing acoustic waves for the Euler equations similar to the traveling wave solution which corresponds to the standing acoustic wave solution.

	It is difficult to identify precisely which limiting solution corresponds to which initial condition. One of the parameters is the value of the initial energy because the energy of the limiting solution cannot exceed the initial energy. Another important parameter is the strength of the resonance between the initial condition and the kernel determined by the shape of the initial condition. The amount of the dispersion is determined by the resonance; in the absence of the resonance the dispersion ceases to be zero and the limiting solution is $ \sigma\equiv 0.$ 

At last, we summarize the main results of our numerical experiments for real as well as the ideal gas background.

In all cases, three types of limiting solutions exist and are attracting in nature.\\
1. The trivial solution $ \sigma\equiv 0 $\\
2. The maximum amplitude traveling waves \eqref{a.2.37} as depicted in Fig.(\ref{fluxb0})\\
3. A family of quasiperiodic solution as represented in Fig.(\ref{fig 6}).

The smooth traveling waves $ |\alpha|\leq 0 $ do not attract solutions and are neutrally stable. A striking feature of the attracting solution is that they all have corners except trivial solution.

\section{Conclusions}\label{con}
In this article, we explored how the real gas effects influence the long-time evolution of a compressible gas modelled by Euler equations of gasdynamics inside a duct. We used an asymptotic method to reduce the Euler equations of gasdynamics in a weakly nonlinear regime with real gas background and reflecting boundary conditions to an integro-differential equation for acoustic components, which is composed of an inviscid burger like nonlinear equation with a linear integral self coupling term and periodic boundary conditions, hence extended the results from ideal gas \cite{MR1719749, MR2716775}  to a van der Waal gas.

We have obtained an exact expression for a family of traveling wave solutions in the case of real gases as well as ideal gas. In each case, the family depends on a parameter 
$ \alpha $ and is composed of smooth wave ($ |\alpha|<1 $) solutions and a maximum amplitude wave ($ \alpha=1 $) solution having a slope discontinuity as depicted in Fig.(\ref{fluxb0}). 

The effects of van der Waals parameter are investigated on the speed of propagation and the amplitude of the traveling waves and are illustrated in Figure (\ref{100}). For a fixed $ \delta=0.4,$ an increase in the van der Waal parameter $ b $ causes the speed of propagation (which is always negative) to decrease. For each value of $ b $, there exists an amplitude range up to a critical value $ 0<A<A_{c}$ and the range of amplitude of the traveling wave decreases with the increase in $ b.$ For different values of van der Waal parameter $ b $  families of traveling wave profiles (for $ \alpha =0.2,\,0.6,\,\text{and}\,1 $) 
are depicted in Figure (\ref{fluxb0}). In each case the traveling waves are smooth for $ |\alpha|<1 $ and have a cusp for $ \alpha=1. $

Evolutionary stages of an arbitrary inital data is displayed in Figs. (\ref{fig131}),(\ref{fig4}),
and (\ref{fig5}), which suggested the weakly dispersive behaviour of the convolution term in evolution equation. In addition to the traveling waves, we have a family of nonbreaking for all time solutions in the real gas background thus extend the ideal gas observations. In spite of having a nontrivial acoustic component, these solutions do not show breaking which leads to shock formation, but may have corner type singularities and shown in Figs. (\ref{fig 6}), (\ref{fig7}). These solutions are different from the traveling wave solutions as unlike traveling wave solutions they change shape during the evolution and forms a globally attracting set for the long time evolution of solution as shown in Figures  (\ref{fig131}),(\ref{fig4}), and (\ref{fig5}).



\end{document}